\documentclass[12pt]{article}
\usepackage{amsfonts}
\usepackage{eepic}
\usepackage{epic}
\usepackage{color}
\usepackage{graphics}
\begin{document}
\newtheorem{Theoreme}{Th\'eor\`eme}[section]
\newtheorem{Theorem}{Theorem}[section]
\newtheorem{Th}{Th\'eor\`eme}[section]
\newtheorem{De}[Th]{D\'efinition}
\newtheorem{Pro}[Th]{Proposition}
\newtheorem{Lemma}[Theorem]{Lemma}
\newtheorem{Proposition}[Theoreme]{Proposition}
\newtheorem{Lemme}[Theoreme]{Lemme}
\newtheorem{Corollaire}[Theoreme]{Corollaire}
\newtheorem{Consequence}[Theoreme]{Cons\'equence}
\newtheorem{Remarque1}[Theoreme]{Remarque}
\newtheorem{Convention}[Theoreme]{{\sc Convention}}
\newtheorem{PP}[Theoreme]{Propri\'et\'es}
\newtheorem{Conclusion}[Theoreme]{Conclusion}
\newtheorem{Ex}[Theoreme]{Exemple}
\newtheorem{Definition}[Theorem]{Definition}
\newtheorem{Remark1}[Theorem]{Remark}
\newtheorem{Not}[Theoreme]{Notation}
\newtheorem{Nota}[Theorem]{Notation}
\newtheorem{Propo}[Theorem]{Proposition}
\newtheorem{exercice1}[Th]{Lemme-Confi\'e au lecteur}
\newtheorem{Corollary}[Theoreme]{Corollary}
\newtheorem{PPtes}[Th]{Propri\'et\'es}
\newtheorem{Defi}[Theorem]{Definition}
\newtheorem{Example1}[Theorem]{Example}
\newenvironment{Proof}{\medbreak{\noindent\bf Proof }}{~{\hskip
3pt$\bullet$\bigbreak}}

\newenvironment{Demonstration}{\medbreak{\noindent\bf D\'emonstration
 }}{~{\hskip 3pt$\bullet$\bigbreak}} 

\newenvironment{Remarque}{\begin{Remarque1}\em}{\end{Remarque1}} 
\newenvironment{Remark}{\begin{Remark1}\em}{\end{Remark1}}
\newenvironment{Exemple}{\begin{Ex}\em}{~{\hskip
3pt$\bullet$}\end{Ex}} 
\newenvironment{exercice}{\begin{exercice1}\em}{\end{exercice1}}
\newenvironment{Notation}{\begin{Not}\em}{\end{Not}}
\newenvironment{Notation1}{\begin{Nota}\em}{\end{Nota}}
\newenvironment{Example}{\begin{Example1}\em}{~{\hskip
3pt$\bullet$}\end{Example1}}

\newenvironment{Remarques}{\begin{Remarque1}\em \ \\* }{\end{Remarque1}}
\newcommand{\Sp}{{\mathbb S}}
\renewcommand{\Re}{{\cal R}}
\renewcommand{\Im}{{\cal F}}
\newcommand{\finpreuve}{~{\hskip 3pt$\bullet$\bigbreak}}
\newcommand{\hp}{\hskip 3pt}
\newcommand{\hph}{\hskip 8pt}
\newcommand{\hphh}{\hskip 15pt}
\newcommand{\vp}{\vskip 3pt}
\newcommand{\vpv}{\vskip 15pt}
\newcommand{\IP}{{\mathbb{IP}}}
\newcommand{\rd}{{\mathbb R}^2}
\newcommand{\R}{{\mathbb R}}
\newcommand{\Hyper}{{\mathbb H}}
\newcommand{\Int}{{\mathbb I}}
\newcommand{\Boule}{{\mathbb B}(0,1)}
\newcommand{\Cantor}{{\mathbb K}}
\newcommand{\K}{{\mathbb K}}
\newcommand{\B}{{\mathbb B}}
\newcommand{\Z}{{\mathbb Z}}
\newcommand{\Nat}{{\mathbb N}}
\newcommand{\N}{{\mathbb N}}
\newcommand{\p}{{\mathbb P}}
\newcommand{\Esp}{{\mathbb E}}
\newcommand{\Complex}{{\mathbb C}}
\newcommand{\Ha}{{\cal H}}
\newcommand{\Harm}{{\bold H}}
\newcommand{\Lcal}{{\cal L}}
\newcommand{\ds}{\displaystyle}
\newcommand{\un}{\bold 1}
\newcommand{\Cone}{C(x,r,\epsilon ,\Phi)}
\newcommand{\Cn}{C(x,2^{-n},\epsilon ,\Phi )}
\newcommand{\Tranche}{W(x,r,\epsilon,\Phi)}
\newcommand{\Wn}{W(x,2^{-n},\epsilon,\Phi)}
\newcommand{\WFn}{W(x,2^{-n},\epsilon,\Phi)\cap F}
\newcommand{\ovec}{\overrightarrow}
\newcommand{\red}{{\bold R}}
\newcommand{\dimH}{\dim_{\Ha}}
\newcommand{\diam}{\mbox{diam}}
\newcommand{\diamit}{\mbox{\em diam}}
\newcommand{\para}{\vskip 2mm}
\newcommand{\cod}{\stackrel{\mbox{\tiny cod}}{\sim}}
\newcommand{\cardit}{\mbox{\em card}}
\newcommand{\card}{\mbox{card}}
\newcommand{\Sphere}{{\mathbb S}_d}
\newcommand{\dist}{\mbox{ dist}}
\newcommand{\distit}{\mbox{\em dist}}
\newcommand{\Tri}{{\cal P}}
\newcommand{\LL}{{\mathcal L}}
\newcommand{\infess}{\mbox{inf\,ess}}
\newcommand{\supess}{\mbox{sup\,ess}}

\definecolor{darkblue}{rgb}{0,0,.5}
\def\u{\underline}
\def\o{\overline}
\def\h{\hskip 3pt}
\def\hh{\hskip 8pt}
\def\hhh{\hskip 15pt}
\def\v{\vskip 8pt}
\def\vv{\vskip 15pt}
\font\courrier=cmr12
\font\grand=cmbxti10
\font\large=cmbx12
\font\largeplus=cmr17
\font\small=cmbx8
\font\nor=cmbxti10
\font\smaller=cmr8
\font\smallo=cmbxti10
\openup 0.3mm
\author{Athanasios BATAKIS and Michel ZINSMEISTER}
\title{On the time schedule of Brownian Flights}
\maketitle

\begin{center}\begin{minipage}{15cm}
\noindent {\bf Abstract:} \em We are interested in the statistics of the duration of Brownian diffusions started at distance $\epsilon$ from the boundary of a given domain and stopped when they hit back this boundary.
\end{minipage}
\end{center}
\section{Introduction} 

The motivation of the following work has its origin in experimental physics. Some long molecules are solvable in a liquid (for instance imogolite in water or DNA in lithium) and the molecules forming the liquid show an intermittent dynamics, alternating diffusion in the bulb and adsorption on the long molecules. For the physicist's point of view, it is very important to have as precise as possible knowledge of the statistics of these brownian flights. 

In \cite{GKLSZ} a connection is established between the statistics of the long flight lengths and the geometry of the long molecules (more precisely their Minkowski dimension).  This connection has been made rigorous in \cite{BLZ}. These two papers concern almost exclusively lengths. How does one check experimentally the results? A very powerful tool for that is relaxation methods in nuclear magnetic resonance (see \cite{CDLPZ}): but this method only allows to compute (the statistics of the) duration of long flights. Some heuristic link between time and length was derived  in \cite{GKLSZ}, \cite{CDLPZ}. The aim of this paper is to make this heuristics rigorous.
 
 \section{Geometric Backgound}
 In the sequel, $\Omega$ will always denote a domain in $\R^d$ with compact boundary.  The crucial tool we need to use is the notion of Whintey cubes. We thus recall the 
 \begin{Proposition}\label{Whitney}  (cf. \cite{Grafakos}, p. 463)
Given any non-empty open proper subset $\Omega$ of $\R^d$, there exists a family of closed dyadic cubes $\{Q_j\}_j$ such that
\begin{itemize}
\item $\bigcup_j Q_j=\Omega$ and the cubes $Q_j$'s have disjoint interiors
\item $\sqrt{d}\ell(Q_j)\le\dist(Q_j,\partial\Omega)\le 4\sqrt{d}\ell(Q_j)$
\item if $Q_j$ and $Q_k$ touch then $\ell(Q_j)\le 4\ell(Q_k)$
\item for a given Whitney cube $Q_j$ there are at most $12^d$ Whitney cubes $Q_k$'s that touch $Q_j$.
\end{itemize}
\end{Proposition}

In this statement, $\ell(Q)$ stands for the side-length of the cube $Q$ and, for $\lambda>0$,  $\lambda Q$ is the cube of the same center and of sidelength $\lambda\ell(Q)$.
For $k\in\Z$, we denote by ${\mathcal Q}_k$, the collection of Whitney cubes $Q_j$ with $\ell(Q_j)=2^k$.
We also recall the definition of the Minkowski sausage: for $r>0$, $$M_r =\{x\in\Omega\;;\; \dist(x,\partial\Omega)\le r\}$$ and
 $$\Gamma_r=\{x\in\Omega\; ;\;\dist(x,\partial\Omega)=r\}$$

We then define  ${\mathcal S}_r$ as the collection of Whitney cubes intersecting $\Gamma_r$. Notice that ${\mathcal S}_r$ is a finite set.

\begin{Definition}
Let $\varepsilon>0$. We will call Brownian flight the random process $F_t, t\ge 0$ consisting in picking at random with equiprobability one of the dyadic Whitney cubes of  ${\mathcal S}_{\varepsilon}$  and starting from the center of the cube a Brownian motion $B_t$ killed once it reaches $\partial\Omega$. We denote by $\tau_{\Omega}= \inf\{t\; ;\; F_t\notin\Omega\}$ the lifetime of this process.
\end{Definition}

We are interested in the asymptotics of $\p(\tau_{\Omega}>t)$ as $t$ grows, but this needs some explanation: 

It is well known that, if $\Omega$ is bounded, this quantity decreases exponentially, as $t\rightarrow\infty$, as $e^{\lambda t}$, where $\lambda$ is the first eigenvalue of the Laplacian. We define 
$$R_{\Omega}=\min\left\{1,\sup_{x\in\Omega} \dist(x,\partial\Omega)\right\}.$$
 Our aim is to evaluate $\p(\tau_{\Omega}>t)$ in the interval $\epsilon^2\le t\le R_{\Omega}^2$, and this independently of $\varepsilon$. In fact, the estimate we are looking for is thus an estimate with respect to $\varepsilon$ rather than for "pure" $t$. 

We study first a simple example, since we will need the partial result anyhow. 
Let $[0,a]$ be a real segment and take $x\in(0,a)$. The probability that brownian motion started at $x$ has not exit the interval $(0,a)$ by time $t$, $\p(\tau_x>t)$, is given by the following equivalent formulas (see \cite{Feller}, pg. 342)
$$\p(\tau_x>t)=\frac{1}{\sqrt{2\pi}}\int_{-x/\sqrt{t}}^{-x/\sqrt{t}}\exp\left({-{\frac12}y^2}\right)dy+ {\sqrt{\frac{2}{\pi}}}\sum_{k=1}^{\infty}(-1)^k\int_{(ka-x)\sqrt{t}}^{(ka+x)\sqrt{t}}\exp\left({-{\frac12}y^2}\right)dy $$
and 
$$\p(\tau_x>t)=\frac{4}{\pi} \sum_{k=1}^{\infty}\exp\left({-\frac{(2n+1)^2\pi^2}{2a^2}t}\right)\sin\frac{(2n+1)\pi x}{a}.$$
By symmetry we can assume $x\le \frac{a}{2}$. If $\frac{a}{\sqrt t}$ is not too big, say $\frac{a}{\sqrt t}<1/2$, we have an easy estimate of $\p(\tau_x>t)$ using the first formula :
$\p(\tau_x>t)\sim \frac{x}{\sqrt t}$.
On the other hand, for $\frac{a}{\sqrt t}$ not too small we get, by the second formula, that 
$\p(\tau_x>t)\le \exp\left({-\frac{\pi^2}{2a^2}t}\right)$. Hence, that there is a change of regime, the decay of  $\p(\tau_x>t)$ with $t$ going from polynomial to exponential and the quantity $\frac{x}{\sqrt t}$ is relevant for small times. 

In a higher dimensional context, let $Q$ be the cube centered at $0$ and of side $r$. By the preceeding remark it follows that, $T$ denoting the  exit time from $Q$ of the Brownian motion starting at $0$, we have
\begin{equation}\label{cubeestimates}
\p(T>t)\le c\left(\frac{r}{\sqrt{t}}\right)^d 
\end{equation}
where $c$ depends only on $d$. Here, $\varepsilon=r\simeq R_{\Omega}$ and the opposite inequality thus holds for  $\varepsilon^2\le t\le R_{\Omega}^2$.

If the Brownian motion is started at distance $\varepsilon$ from the boundary then the exit time is essentially the exit time from a half space and we thus get 
\begin{equation}\label{cubetimetable}
\p(T>t)\le c\left(\frac{\varepsilon}{\sqrt{t}}\right)= c\left(\frac{\varepsilon}{\sqrt{t}}\right)^{d-1+2-d},
\end{equation}
and, as we will see, the opposite inequality is valid for $t$ not too big.

Our goal is to extend (\ref{cubetimetable}) to general domains with rough boundary. In order to describe the domain of validity of our result let us recall a few definitions. Let $K$ be a compact subset of $\R^d$. For $j\ge 0$ let $N_j$ be the number of dyadic cubes of the $j$-th generation (i.e of size $2^{-j}$) that intersect $K$.

\begin{Definition}
The Minkowski dimension of $K$ is 
$$d_M(K)=\limsup_{j\to\infty}\frac{\log_2(N_j)}{j}$$
\end{Definition}

Returning to our situation, we can define similarly the Whitney dimension of $\partial \Omega$ as 
\begin{equation}\label{Whitneydimension}
d_W=d_W(\partial \Omega)=\limsup_{j\to\infty}\frac{\log_2(W_j)}{j},
\end{equation}
where $W_j$ is the number of elements of ${\mathcal Q}_j$.

Under very mild conditions (see   \cite{Bishop3}, \cite{JK}, \cite{BLZ}) these two dimensions coincide.
If the boundary of $\Omega$ has some self similarity we can moreover say that there is a constant $c>0$ such that 
\begin{equation}\label{hypothesis}
\frac1c{\varepsilon}^{d_M}\le\#{\mathcal S}_{\varepsilon}\le c {\varepsilon}^{d_M},
\end{equation}
for all $\varepsilon\le R_{\Omega}$, where $d_M=d_M(\partial \Omega)$. 

In our main theorem we will assume that our domain $\Omega$ satisfies (\ref{hypothesis}).

We also suppose that the domain $\Omega$ satisfies so-called $\Delta$-regularity condition (see also \cite {JW}, \cite{Ancona5}, \cite{hkm}): there exists $L>0$ such that for all $x\in\Omega$, if $d_x=\dist(x,\partial\Omega)<R_{\Omega}$ then 
\begin{equation}\label{capacitycondition}
\ds\omega_{\B(x,2d_x)\cap\Omega}^x\left(\partial\Omega\right)\ge L,
\end{equation} 
where $\ds\omega_{\B(x,2d_x)\cap\Omega}^x$ is the distribution law of the hitting point of Brownian motion starting at x and killed when reaching the boundary of  $\B(x,2d_x)\cap\Omega)$. 
This is a very mild condition (satisfied, for instance, by all domains in $\R^2$ with non-trivial connected boundary) that appears frequently in related literature in various forms (for instance ``uniform capacity condition'' or Hardy inequality).

\section{Time and length estimates for Brownian flights}
We can now state the main result. Let $\Omega$ be a bounded domain in $\R^d$ satisfying (\ref{hypothesis}) and (\ref{capacitycondition}). If $\tau_{\Omega}$ denotes the life-time of a Brownian flight $F_t$ with parameter $\varepsilon$ we have 

\begin{Theorem}\label{Th1}
There exists $c>0$ depending only on constants in  (\ref{hypothesis}), (\ref{capacitycondition}) (and in particular not on $\varepsilon$) such that  \begin{equation}\label{Th1.1}
\frac1c\left(\frac{\varepsilon}{\sqrt t}\right)^{d_{M}+2-d}\le\p(\tau_{\Omega}>t)
\end{equation}
and 
\begin{equation}\label{Th1.2}
\p(\tau_{\Omega}>t)\le c  \left(\frac{\varepsilon}{\sqrt t}\right)^{d_{M}+2-d}\left|\log\left(\frac{\varepsilon}{\sqrt t}\right)\right|^{2d},
\end{equation}
for all $\varepsilon^2<t<R_{\Omega}^2$.
\end{Theorem}

This theorem has a ``cousin'' theorem, which was proved in \cite{BLZ}.

\begin{Theorem}\label{main2}
Let $\varepsilon<r<R_{\Omega}$. The probability that the hitting point of $F$ is at distance greater than $r$ from the starting point $x$ is comparable to 
\begin{equation}\label{main2.1}
\left(\frac{\#{\mathcal S}_r}{\#{\mathcal S}_{\varepsilon}}\right)^{d_M}\left(\frac{r}{\varepsilon}\right)^{d-2}
\end{equation}
\end{Theorem}

Notice that we do not assume ($\ref{hypothesis}$) for this theorem. If we do, we have
\begin{equation}\label{main2.2}
\left(\frac{\#{\mathcal S}_r}{\#{\mathcal S}_{\varepsilon}}\right)^{d_M}\left(\frac{r}{\varepsilon}\right)^{d-2}\sim \left(\frac{r}{\varepsilon}\right)^{d_M-(d-2)}
\end{equation}

Notice that the quantity on the left of (\ref{main2.2}) is the same as the one in (\ref{Th1.1}) where we have replaced $r$ by $\sqrt{t}$ which is coherent with standard behaviour of Brownian motion.

\section{Proof of the theorem \ref{Th1}}

For $s>0$ we denote by $\beta_s$ the total time  spent by Brownian flight $F_t$  in the Minkowski sausage $\{x\in\Omega\;;\; \dist(x,\partial\Omega)\le s\}$ and 
$\delta\beta_s=\beta_s-\beta_{s/2}$.

We define analogue quantities more adapted to the Whitney decomposition ; namely, $\delta \tilde\beta_{2^k}$ will denote the time spent by $F_t$ inside ${\tilde M}_k=\bigcup\{Q\;;\;Q\in{\mathcal Q}_k\}$.

If $Q$ is a Whintney cube we define the "vicinity" of $Q$ as $$\tilde Q=Q\cup\bigcup\{Q'\;;\;Q'\in{\mathcal E}\},$$ 
where $${\mathcal E}=\{Q'\in\bigcup_k{\mathcal Q}_k\;;\; \lambda Q\cap Q'\not= \emptyset \mbox{ and } \lambda Q'\cap Q\not=\emptyset\},$$
and $\lambda= 8\sqrt{d}$ satisfy that for all Whitney cubes 
$$\lambda Q\supset \B(x_Q,2\dist(x_Q,\partial\Omega)),$$ $x_Q$ being the center of $Q$.

We may now start the proof and we begin with 

\subsection{The upper bound}

We separate the event $\{\tau_{\Omega}>t\}$ by the partition $\{\tau_{\Omega}=\beta_{\sqrt{t}}\}$ and $\{\tau_{\Omega}>\beta_{\sqrt{t}}\}$, that is whether the process goes or does not go at distance $\sqrt{t}$ from the boundary. From theorem \ref{main2} we get
$$\p(\tau_{\Omega}>t\mbox{ and }\tau_{\Omega}>\beta_{\sqrt{t}})\le \p(\tau_{\Omega}>\beta_{\sqrt{t}})\le c \left(\frac{\varepsilon}{\sqrt{t}}\right)^{d_M+2-d}.$$

In order to estimate the term $\p(\tau_{\Omega}=\beta_{\sqrt{t}}>t)$ we begin by writing $$\ds\tau_{\Omega}=\sum_{k=-\infty}^{\log_2{\sqrt{t}}}\delta\beta_{2^k},$$
and thus, putting $k_0=\log_2{\epsilon}$ and $k_1=\log_2{\sqrt{t}}$,
\begin{equation}\label{sum}
\p\left( \sum_{k=-\infty}^{\log_2{\sqrt{t}}}\delta\beta_{2^k} >t \right)\le \sum_{k=k_0+1}^{k_1}\p\left(\delta\beta_{2^k}>\frac{t}{2(k_1-k_0)} \right)+\p\left(\beta_{2^{k_0}}>\frac{t}{2}\right)
\end{equation}

We now invoke the following lemma whose proof is post poned to the next section.

\begin{Lemma}\label{Technical1} There exists a number $k^*$ depending only on $d$ and constants $C$ , $0<p<1$ depending only on $d$ and $L$ such that for all $t>0$, $k\in\Z$ and $N\in\N$  we have 
\begin{eqnarray*}
&&\p(\delta\beta_{2^k}>t)\le Cp^N 
+\\
&&\p\left(\exists Q \in{\mathcal Q}_k\cup...\cup {\mathcal Q}_{k-k^*}; \exists\; 0<s_1<s_2<t \mbox{ with }F_{[s_1,s_2]}\subset\tilde Q \mbox{ and } s_2-s_1>t/N\right)
\end{eqnarray*}
\end{Lemma}

Using this lemma, we get
 \begin{eqnarray*}
&&\p\left(\delta\beta_{2^j}>\frac{t}{2(k_1-k_0)}\right)\le cp^N+\\
&&\p\left(\exists Q\in{\mathcal Q}_{j}\cup...\cup{\mathcal Q}_{j-k^*}\; \exists s_1<s_2<\tau_{\Omega}\;;\;F_{[s_1,s_2]}\subset\tilde Q\;,\; s_2-s_1>\frac{t}{2(k_1-k_0)N}\right)
\end{eqnarray*}

By (\ref{cubeestimates}), (\ref{main2.2}) and the strong Markov property of Brownian motion we then have :

 \begin{eqnarray*}
&&\p\left(\delta\beta_{2^j}>\frac{t}{2(k_1-k_0)}\right)\le cp^N+\\
&&\p\left(\exists Q\in{\mathcal Q}_{j}\cup...\cup{\mathcal Q}_{j-k^*}\; \exists s_1\tau_{\Omega}\;;\;F_{s_1}\in\tilde Q\right)\times \\
&&\p\left(\exists s_2>s_1\;,\; F_{[s_1,s_2]}\subset\tilde Q\mbox{ and }s_2-s_1>\frac{t}{2(k_1-k_0)N}\right)\le \\
&& cp^N+c\left(\frac{2^j}{\sqrt{\frac{t}{2(k_1-k_0)N}}}\right)^d\left(\frac{\varepsilon}{2^j}\right)^{d_M+2-d}
\end{eqnarray*}

Suming up the first term in (\ref{sum}) we get
 \begin{eqnarray*}
&&\sum_{k=k_0+1}^{k_1}\p\left(\delta\beta_{2^k}>\frac{t}{2(k_1-k_0)} \right)\le c(k_1-k_0)p^N+
c\sum_{k=k_0+1}^{k_1}\left(\frac{2^k}{\sqrt{\frac{t}{2(k_1-k_0)N}}}\right)^d \left(\frac{\varepsilon}{2^k}\right)^{d_M+2-d}\\
&&\le c(k_1-k_0)p^N+c(2(k_1-k_0)N)^d\varepsilon^{d_M+d-2}\left(\frac{1}{\sqrt{t}}\right)^d\sum_{k=k_0+1}^{k_1} 2^{2d-d_M+2}\\
&&\le  c(k_1-k_0)p^N+c(2(k_1-k_0)N)^d\varepsilon^{d_M+d-2}\left(\frac{1}{\sqrt{t}}\right)^d 2^{k_1(2d-d_M+2)}\\
&&\le  cp^N\log_2\left(\frac{\sqrt{t}}{\varepsilon}\right)+c\left(\log_2\left(\frac{\sqrt{t}}{\varepsilon}\right)N\right)^d\left(\frac{\varepsilon}{\sqrt{t}}\right)^{d_M+2-d}
\end{eqnarray*}

Take $N\simeq (d_M+2-d) \log_p\left(\frac{\varepsilon}{\sqrt{t}}\right)$ to obtain
\begin{equation}\label{firstterm}
\sum_{k=k_0+1}^{k_1}\p\left(\delta\beta_{2^k}>\frac{t}{2(k_1-k_0)} \right)\le c \left(\log\left(\frac{\varepsilon}{\sqrt{t}}\right)\right)^{2d} \left(\frac{\varepsilon}{\sqrt{t}}\right)^{d_M+2-d}
\end{equation}

To bound the second term $\p\left(\beta_{2^{k_0}}>\frac{t}{2}\right)$ of the sum (\ref{sum}) we need a lemma of the same nature as lemma \ref{Technical1}. The proof of this lemma is also post poned to the next section.
\begin{Lemma}\label{Technical4}
Let ${\mathcal R}_{k_0}$ be the collection of all dyadic cubes of sidelength $2^{k_0}$ intersecting $\partial\Omega$. There exist constants  $C$ , $0<p<1$ depending only on $d$ and $L$ such that for all $t>0$, $k_0\in\Z$ and $N\in\N$  
\begin{eqnarray*}
&&\p(\delta\beta_{2^{k_0}}>t)\le Cp^N 
+\\
&&\p\left(\exists Q \in{\mathcal R}_{k_0}; \exists\; 0<s_1<s_2<t \mbox{ with }F_{[s_1,s_2]}\subset NQ \mbox{ and } s_2-s_1>t/N\right)
\end{eqnarray*}
\end{Lemma}

We therefore deduce 
$$\p(\delta\beta_{2^{k_0}}>\frac{t}{2})\le cp^N+\left(\frac{2^{k_0}N}{\sqrt{t}}\right)^d.$$
We minimize on $N\simeq\log_p(\varepsilon)$ to get
$$\p(\delta\beta_{2^{k_0}}>\frac{t}{2})\le c\left(\log(\varepsilon)\frac{\varepsilon}{\sqrt{t}}\right)^d.$$

Combining this last inequality and (\ref{firstterm}) we get the upper bound, since $d_M+d-2\le d$.

\subsection{The Lower Bound}\label{lb}
Following the same reasoning for ${k_1}=[\log_2\sqrt{t}]+1$ we get 
$$\p(\tau_{\Omega}>t)\ge \p(\exists s_1>0 \mbox{ s.t. }F_{s_1}\in\bigcup_{Q\in{\mathcal Q}_{k_1}}Q\mbox{ and }\exists {s_2>s_1+t} \mbox{ s.t. } F_{[s_1,s_2]}\in\bigcup_{Q\in{\mathcal Q}_{k_1}}2 Q)$$
Using strong Markov property the this probability can be written as the product of 
$\p(\exists s_1>0 \mbox{ s.t. }F_{s_1}\in{\mathcal Q}_{k_1} \mbox)$ with $\p(\exists {s_2>s_1+t} \mbox{ s.t. } F_{[s_1,s_2]}\in\bigcup_{Q\in{\mathcal Q}_{k_1}}2 Q)$. The second term of the product  is greater than the probability that Brownian motion exits a cube of size $2^{k_1+1}\simeq\sqrt{t}$ at time greater that $t$ which is bounded below by a positive constant depending only on $d$. The first one is simply the probability that Brownian flight gets to ${\mathcal Q}_{k_1}$ which is equivalent to $\left(\frac{\epsilon}{\sqrt{t}}\right)^{d_{M}+2-d}$ and the proof is complete.

\section{Proofs of lemmas}
Let un first deal with lemma  \ref{Technical1}. The proof of  \ref{Technical4} is quite similar and will hence be abridged.
\subsection{Proof of lemma \ref{Technical1}}
We need the following
\begin{Lemma}\label{key} Under the $\Delta$-regularity hypothesis
the probability that BM touches more than $N$ Whitney cubes of a given size decreases as $Cp^N$, with $0<p<1$, $C$ a positive constant.
\end{Lemma}

 The proof of the lemma relies on an annuli reasoning.
 \begin{Proof}
 Le $(B_t)_{t>0}$ be Brownian motion started at any point $x\in\Omega$ and choose $k\in\Z$. Choose any $Q\in{\mathcal Q}_k$ and let $\lambda Q$ be the cube of the same center but $\lambda$ times the side-length $\ell(Q)$  of $Q$. By the definition of Whitney cubes, there is a $\lambda= 8\sqrt{d}$ depending only on $d$  such that $$\frac{\lambda}{2} \ell(Q)\le\dist(Q,\partial\Omega)\le 2\lambda \ell(Q).$$ 

Suppose that there exists $t_0>0$ such that $B_{t_0}\in Q$.
By the $\Delta$-regularity condition (\ref{capacitycondition}), the probability that there exists $t_1>t_0$ with $B_{[t_0,t_1]}\subset\Omega$ and $B_{t_1}\notin\lambda Q$ is bounded above by $p<1$ depending only on $L,\lambda$:
\begin{equation}\label{one}
\p\left(\exists t_1>t_0\; ; \; B_{[t_0,t_1]}\subset\Omega \mbox{ and }B_{t_1}\notin\lambda Q\Big|\exists t_0>0\;;\; B_{t_0}\in Q\right)<p
\end{equation}
 
On the other hand, the number of Whitney cubes of ${\mathcal Q}_k$ lying inside $\lambda Q$ is bounded by a constant $c_1=c_1(d)$. The probability that there exists a Whitney cube $Q_1\in {\mathcal Q}_k$ outside $\lambda Q$ that is visited by Brownian motion is hence bounded above by $p<1$. 

We study probability that there exist Whitney cubes $Q_1,...,Q_m\in{\mathcal Q}_k$ such that $Q_1\cap\lambda Q=Q_2\cap\lambda Q_1=...=Q_m\cap\lambda Q_{m-1}=\emptyset$ all visited by Brownian motion.
It is sufficient to prove that this probability decays exponentially with $m$.

By the strong Markov property the probability that there exists $t_m>t_{m-1}>...>t_0$ such that $B_{t_0}\in Q$ , $B_{t_1}\in Q_1$ ... , $B_{t_m}\in Q_m$ is given by 
\begin{eqnarray*}
&&\p\left(\exists\; t_m>t_{m-1}>...>t_0\; ;\mbox{ and }Q,...Q_m \mbox{ as above such that } B_{t_m}\in Q_m,...,B_{t_0}\in Q\right)\\
&&=\p\left(\exists\; t_m>t_{m-1}\; ;\;B_{t_m}\in Q_m|\exists\; t_{m-1}>...>t_0\; B_{t_{m-1}}\in Q_{m-1},...,B_{t_0}\in Q\right)\cdot\\
&&\p\left( \exists\; t_{m-1}>...>t_0\; B_{t_{m-1}}\in Q_{m-1},...,B_{t_0}\in Q\right)\\
&&=\p\left(\exists\; t_m>t_{m-1}\; ;\;B_{t_m}\in Q_m \mbox{ with }Q_m\cap\lambda Q_{m-1}=\emptyset|\exists\; t_{m-1}\; ;B_{t_{m-1}}\in Q_{m-1}\right)\cdot\\
&&\p\left( \exists\; t_{m-1}>...>t_0\; B_{t_{m-1}}\in Q_{m-1},...,B_{t_0}\in Q\right)
\end{eqnarray*}
Now, by (\ref{one}),  
$$\p\left(\exists\; t_m>t_{m-1}\; ;\;B_{t_m}\in Q_m \mbox{ with }Q_m\cap\lambda Q_{m-1}=\emptyset|\exists\; t_{m-1}\; ;B_{t_{m-1}}\in Q_{m-1}\right)<p.$$ 
By induction we get that 
$$\p\left(\exists\; t_m>t_{m-1}>...>t_0\; ;\mbox{ and }Q,...Q_m \mbox{ as above such that } B_{t_m}\in Q_m,...,B_{t_0}\in Q\right)<p^m$$ and hence the lemma.
\end{Proof}

Recall that for a given dyadic Whitney cube $Q $ we have defined the vicinity $\tilde Q$ of $Q$ as the union of all Whitney cubes $Q'$ verifying 
$$Q'\cap\lambda Q\not=\emptyset\mbox{ and }Q\cap\lambda Q'\not=\emptyset.$$ 
We can easily check that there are less than $(100\sqrt{d})^d$ such cubes $Q'$ of size at most $\ell(Q)/12$  (the constants are not optimal). 
We say that the $k$-level layers are visited at least $n$ times if there exist $t_0<s_1<t_1<...<s_n<t_n$ satisfying 
$$B_{t_j}\in\bigcup_{Q\in{\mathcal S}_{2^k}}Q\mbox{ and }B_{s_j}\notin\bigcup_{Q\in{\mathcal S}_{{2^k}}}\tilde Q,$$ for all $j=1,...,n$. For any $k\in\Z$ note 
$$\nu_k=\sup\{n\in\N\; ;\; \mbox{the $k$-level layers are visited at least $n$ times}\}$$
\begin{Lemma}\label{technical0}
There exists $0<p<1$ and a positive constant $C$ such that, given $k\in\Z$, for all $n\in\N$
$$\p(\nu_k>n)\le p^n \p\left(\exists t_0>0\mbox{ and } Q\in{\mathcal S}_{2^k}\; ; \; B_{t_0}\in Q\right).$$
\end{Lemma}
\begin{Proof}
The arguments as similar as in lemma \ref{key}. We only need to prove that $\p(\nu_k>1)<p$ and apply strong Markov property. 
We have 
\begin{eqnarray*}
\p(\nu_k>1)& \le &\p\left(\exists 0<t_0<s_1<t_1\; ,\; Q\in{\mathcal S}_{2^k} \; ; \; B_{t_0}\in Q\; , \; B_{s_1}\notin \bigcup_{Q\in{\mathcal S}_{{2^k}}}\tilde Q\; , \; B_{t_1}\in\bigcup_{Q\in{\mathcal S}_{2^k}}Q\right)\\
&=& \p\left(\exists t_1>s_1>t_0\; ; \;  B_{s_1}\notin \bigcup_{Q\in{\mathcal S}_{{2^k}}}\tilde Q\; , \; B_{t_1}\in\bigcup_{Q\in{\mathcal S}_{2^k}}Q\Big| \exists t_0>0 \; ; \; B_{t_0}\in Q\in{\mathcal S}_{2^k}\right)\\
&\times &\p\left(\exists t_0>0\mbox{ and } Q\in{\mathcal S}_{2^k}\; ; \; B_{t_0}\in Q\right)
\end{eqnarray*}
To abbreviate formulas we note $ \p_c(.)=\p\left(.| \exists t_0>0 \; ; \; B_{t_0}\in Q\in{\mathcal S}_{2^k}\right)$. With this notation,
\begin{eqnarray*}
&&\p_c\left(\exists t_1>s_1>t_0\; ; \;  B_{s_1}\notin \bigcup_{Q\in{\mathcal S}_{{2^k}}}\tilde Q\; , \; B_{t_1}\in\bigcup_{Q\in{\mathcal S}_{2^k}}Q)\right)=\\
&&\p_c\left(\exists t_1>s_1\; ; \; B_{t_1}\in\bigcup_{Q\in{\mathcal S}_{2^k}}Q\; \Big| A\right)\p_c(A)+\p_c\left(\exists t_1>s_1\; ; \; B_{t_1}\in\bigcup_{Q\in{\mathcal S}_{2^k}}Q\; \Big| B\right)\p_c(B)
\end{eqnarray*}
where $$A=\left\{\exists s_1>t_0\; ; \;  B_{s_1}\notin \bigcup_{Q\in{\mathcal S}_{{2^k}}}\lambda Q\right\}\mbox{ and }$$
$$B=\left\{\exists s_1>t_0\; ; \;  B_{s_1}\in Q'\ , \,\lambda Q'\cap \bigcup_{Q\in{\mathcal S}_{{2^k}}}Q=\emptyset\;, \;B_{s_1}\in \bigcup_{Q\in{\mathcal S}_{{2^k}}}\lambda Q\right\}$$
form a partition of the event $\left\{\exists s_1>t_0\; ; \;  B_{s_1}\notin \bigcup_{Q\in{\mathcal S}_{{2^k}}}\tilde Q\right\}$.


By (\ref{one}), $\p_c(A)\le p$. Similarly, by the strong Markov property of Brownian motion and by (\ref{one}),
$$\p_c\left(\exists t_1>s_1\; ; \; B_{t_1}\in\bigcup_{Q\in{\mathcal S}_{2^k}}Q\; \Big| B\right)= \p\left(\exists t_1>s_1\; ; \; B_{t_1}\in\bigcup_{Q\in{\mathcal S}_{2^k}}Q\; \Big| B_{s_1}\in Q' \right)\le p.$$
We deduce that  
$$\p_c\left(\exists t_1>s_1>t_0\; ; \;  B_{s_1}\notin \bigcup_{Q\in{\mathcal S}_{{2^k}}}\tilde Q\; , \; B_{t_1}\in\bigcup_{Q\in{\mathcal S}_{2^k}}Q)\right)\le \p_c(A)+p(1-\p_c(A))$$

The function $t\mapsto t+p(1-t)$ being increasing on $[0,p]$  we get 
$$\p_c\left(\exists t_1>s_1>t_0\; ; \;  B_{s_1}\notin \bigcup_{Q\in{\mathcal S}_{{2^k}}}\tilde Q\; , \; B_{t_1}\in\bigcup_{Q\in{\mathcal S}_{2^k}}Q)\right)\le 2p-p^2<1$$
and the lemma is proven.
\end{Proof}

Remark that by definition of dyadic Whitney cubes there exist $k^*$ depending only on the dimension of the space ($k^*=[\log_2(8\sqrt{d})]+3$ will do) such that for all $k\in\Z$,  
\begin{equation}\label{k*}
\{x\in\Omega\;;\;2^{k-1}\le \dist(x,\partial\Omega)\le 2^k\}\subset\bigcup_{j=k-k^*}^k\bigcup_{Q\in{\mathcal Q}_j}Q.
\end{equation}


\begin{Proof}{\bf of lemma \ref{Technical1}}
We clearly have 
\begin{equation}\label{partition}
\p(\delta\beta_{2^k}>t)\le \p(\delta\beta_{2^k}>t\; ,\;\nu_k+...+\nu_{k-k^*}>N)+\p(\delta\beta_{2^k}>t\; ,\;\nu_k+ ...+\nu_{k-k^*}\le N)
\end{equation}
Given $t>0$, $k\in\Z$,  by lemma \ref{technical0} we have, for all $N$,
$$ \p(\delta\beta_{2^k}>t\; ,\;\nu_k+...+\nu_{k-k^*}>N)< k^*p^{\frac{1}{k^*}N} \p\left(\exists t_0>0\mbox{ and }Q\in{\mathcal S}_{2^{k-k^*}}\; ; \; B_{t_0}\in Q\right)<c\tilde p^N$$
Let us estimate the second term of the sum (\ref{partition}). 
By (\ref{k*}) and the definition of $\nu_k$ we get
 \begin{eqnarray*}
 &&\p(\delta\beta_{2^k}>t\; ,\;\nu_k+...+\nu_{k-k^*}\le N )\\
&&\le\p(\exists l\le N\; , \;Q_1,...,Q_l\in\bigcup_{j={k-k^*}}^k\mathcal{Q}_j\;;\;Q_s\cap\tilde Q_{s-1}=\emptyset\;,\;\forall s=2,...,l\mbox { and } \\
&&\exists  t_1<s_1\le t_2<s_2<...\le t_l<s_l\; ; \;B_{[t_i,s_i]}\subset \tilde Q_s\;\forall s=1,...,l\mbox{ and }\sum_{i=1}^ls_i-t_i>t).
\end{eqnarray*}
Since $l\le N$ we get that 
$\p(\delta\beta_{2^k}>t\; ,\;\nu_k+...+\nu_{k-k^*}\le N ) $ is bounded above by 
$$\p\left(\exists Q \in{\mathcal Q}_k\cup...\cup {\mathcal Q}_{k-k^*}; \exists\; 0<t_{i_0}<s_{i_0}\mbox{ with }B_{[t_{i_0},s_{i_0}]}\subset\tilde Q \mbox{ and } s_{i_0}-t_{i_0}>t/N\right),$$
which completes the proof.
\end{Proof}
\subsection{Proof of lemma \ref{Technical4}}
The ideas are the same but we will work with cubes touching the boundary instead of Whitney cubes. 
\begin{Lemma}\label{key1} Under the $\Delta$-regularity hypothesis, the probability that BM started at distance $\varepsilon<r$ from the boundary gets at distance greater than $R$ from the starting point without leaving the Minkowski sausage $M_r=\{x\in\Omega\;;\; \dist(x,\partial\Omega)\le r\}$ is bounded above by $cp^{R/r}$ where $c>0$ and $0<p<1$ are constants (depending only on $L$ of the $\Delta$-regularity hypothesis and on $d$).
\end{Lemma}
The proof, similar to the one of lemma \ref{key}, is therefore abridged.
\begin{Proof}
Let $x\in M_r$ and consider the annuli centered at $x$ of inner radii $4\ell r$ and outer radii $4(\ell+1)r$ where $\ell=0,...,\left[\frac{R}{4r}\right]$. Brownian motion started at $x$ and moving at distance $R$ from $x$ before exiting $\Omega$ must go through all these annuli. 
The probability of going through such an annulus while staying at distance at most $r$ from the boundary is bounded by a $p_0\in (0,1)$ by the $\Delta$-regularity hypothesis. To see this take any point $y$ in the middle of the annulus (i.e. at distance $\frac{6\ell+6}r$ from $x$) and consider the ball of center $y$ and radius $2r$. If $d_y<r$, the probability to exit the ball without  touching $\partial\Omega$ is uniformly bounded away from $1$ by the same hypothesis. This probability being greater than the probability of going through the annulus  we have the statement.
By the independence of the ``crossing annuli" events we get that the probability that Brownian motion  goes through all the annulli is smaller that $cp_0^{\left[\frac{R}{4r}\right]}\sim c\tilde p^{\left[\frac{R}{r}\right]}$.
\end{Proof}

We say that the Minkowki sausage $M_r$ is visited by the Brownian motion at least $k$ times if there exist $t_0<s_1<t_1<...<s_n<t_n<\tau_{\Omega}$ satisfying $B_{t_i}\in M_r$ for all $i=0,...,n$ and $B_{s_i}\notin M_{4r}$.

In a similar way with $\nu_k$ we define $\xi_k$ as the 
$$\xi_{r}=\sup\{n\in\N\;;\ M_r \mbox{ has been visited at least } k \mbox{ times}\}.$$
\begin{Lemma}\label{technical1}
There exists $0<p<1$ depening only on the $\Delta$-regularity's $L$ such that, given $r>0$, for all $n\in\N$
$$\p(\xi_{r}>n)\le p^n.$$
\end{Lemma}
\begin{Proof}
The proof of this lemma is a straightforward  application of the $\Delta$-regularity condition.
It suffices to show that there exists $0<p<1$ such that
$$\p(\left(\exists t_1>s_1> t_0\; ; \; B_{t_0}\in M_r\, ,\,B_{s_1}\notin M_{4r}\,,\, B_{t_1}\in M_r \right)<p,$$
and then apply  the Markov property. 
Remark that, the probability 
$$\p(\left(\exists t_1>s_1> t_0\; ; \; B_{t_0}\in M_r\, ,\,B_{s_1}\notin M_{4r}\right)$$ is smaller than the probability that brownian motion started at $t_0$ exits a ball of radius $2r>2\dist(B_{t_0},\partial\Omega)$ without hitting $\partial\Omega$. By the $\Delta$-regularity this last probability is bounded by a constant $p<1$.
\end{Proof}
\begin{Proof}{\bf of lemma \ref{Technical4}}
As before we get :
\begin{equation}\label{partitionbis}
\p(\delta\beta_{2^{k_0}}>t)=\p(\delta\beta_{2^{k_0}}>t\; ,\;\xi_{2^{k_0}}>N)+\p(\delta\beta_{2^{k_0}}>t\; ,\;\xi_{2^{k_0}}\le N)
\end{equation}
By lemma \ref{technical1} we get $\p(\delta\beta_{2^{k_0}}>t\; ,\;\xi_{2^{k_0}}>N)\le p^N$, where $0<p<1$, for all $N\in\N$.
Let us now deal with the second term of the sum.
$$\p(\delta\beta_{2^{k_0}}>t\; ,\;\xi_{2^{k_0}}\le N)\le \p(\exists s_1<s_2<\tau_{\Omega}\;;\; s_2-s_1>t/N\,,\,B_{[s_1,s_2]}\subset M_{4r})$$
Using lemma \ref{key1} we get that, for $R>4r$,  this probability is bounded by
$$cp^\frac{R}{r}+  \p\left(\exists s_1<s_2<\tau_{\Omega}\;;\; s_2-s_1>t/N\,,\,B_{[s_1,s_2]}\subset M_{4r}\cap\B(B_{s_1},R)\right),$$ and the statement of the lemma  \ref{Technical4} follows on taking $N=\left[\frac{R}{r}\right]$
\end{Proof}
\section{Further Comments}
We should point out that hypothesis  (\ref{hypothesis})  in theorem \ref{Th1} can be dropped; in this case, the same reasoning as in subsection \ref{lb} gives the lower bound $$\left(\frac{\#{\mathcal S}_{\sqrt t}}{\#{\mathcal S}_{\varepsilon}}\right)^{d_M}\left(\frac{{\sqrt t}}{\varepsilon}\right)^{d-2}.$$ The best upper bound is less evident; nevertheless a slight improvement of the above proofs gives 
$$\p(\tau_{\Omega}>t)\le c\left(\log\frac{\sqrt t}{\epsilon}\right)^{cd}\sup_N\left(2^{-Nd}+\sum_{k=0}^N\left(\frac{\#{\mathcal S}_{\sqrt t/2^k}}{\#{\mathcal S}_{\varepsilon}}\right)^{d_M}\left(\frac{2^{-k}{\sqrt t}}{\varepsilon}\right)^{d-2}\right).$$

\bibliographystyle{alpha}
\bibliography{biblio}

\end{document}